\documentclass[10pt]{article}
\usepackage{amsmath}
\usepackage{amssymb}
\usepackage{amsthm}
\usepackage{graphicx} 
\usepackage{appendix}

\newtheorem{theorem}{Theorem}[section]
\newtheorem{lemma}[theorem]{Lemma}

\newtheorem{definition}[theorem]{Definition}

\newcommand{\Sup}[1]{\underset{#1}{sup\,}}

\newcommand{\R}{{\bf R}}

\newcommand{\Ome}{\Omega}
\newcommand{\lap}{{\Delta}}

\newcommand{\bra}[1]{\left( #1 \right)}
\newcommand{\Quad}[1]{{Q \bra{#1}}}

\newcommand{\norm}[1]{\Vert #1 \Vert}

\newcommand{\Schrodinger}{Schr\"odinger }

\newcommand{\beq}{\begin{equation}}

\newcommand{\mlap}[1]{\norm{\bra{-\lap}^{#1}f}}

\newcommand{\kik}[1]{C^{\infty}_{c}\bra{#1}}
\newcommand{\ip}[2]{\langle #1,#2\rangle}

\newcommand{\sx}{\langle x \rangle}

\newcommand{\sobz}{W^{m,2}_0}

\begin{document}
\title{\bf Pointwise Lower Bounds on the Heat Kernels of Higher Order \Schrodinger Operators}
\author{Narinder S Claire}
\date{}
\maketitle
\begin {abstract}

We apply Davies' method for obtaining pointwise lower bounds on the heat kernels of higher-order differential operators to obtain 
pointwise lower bounds in the presence of a polynomialy bounded potential. 
\\ \\
{\bf AMS Subject Classification :   35K41 }  
\end{abstract}

\section{Introduction}
There has been extensive work on the theory for heat kernels of second order elliptic operators but the
situation for higher order ellipitic operators is considerably different. Higher order heat kernels are not necessarily 
positive on the off-diagonal. Despite the limitations 
of the approaches often used in the case of second order operators, Davies\cite{Da4} provided a powerfull tool
obtaining the necessary bounds and established a lower bound for the heat
kernel of  higher order elliptic operators with bounded measurable coefficients on $L^2\bra{\R^N}$ where  $2m>N$. 
Robinson and Ter Elst\cite{Rob} obtained bounds  for operators on a
Lie Group . Their method however required the conservation of probability.

In this exposition we show that Davies' theory can also be exploited to obtain pointwise lower bounds 
for heat kernels of the Schr\"odinger operator with a polynomialy bounded potential. \\
We define 
$$\sx:=\sqrt{\bra{|x|^2+\rho}} \quad \text{ where } \rho = \max\{1,\bra{\frac{2m}{N}-1}^{\frac{1}{m}}\}$$
We will assume that $H_0$ is an elliptic differential operator of order $2m$ on $L^2\bra{\R^N}$ with constant of ellipticity $c_1$ 
and that for the potential $V$ there are positive constants $c_2$ and $\gamma$ such that $$ V\bra{x} \leq c_2\sx^\gamma$$ 

If $Q_0$ is the quadratic form of $H_0$  defined with domain equal to the Sobolev space 
$\sobz \bra{\R}$ where 
\begin{equation}\label{H1}
c_{1}^{-1}\mlap{\frac{m}{2}}^2_2\leq Q_0\bra{f}\leq c_{1}\mlap{\frac{m}{2}}^2_2
\end{equation}
Then we define $Q$, the quadratic form associated with $H:=H_0+V$  
$$Q\bra{f} := Q\bra{f} + \norm{V^{\frac{1}{2}}f}_2^2$$

Our main result is the following theorem:
\begin{theorem}\label{mian}
If $k\bra{t,x,y}$, the heat kernel generated by $H$ is continuous in all three variables and
is in  $L^2\bra{\R^N}$ for all $x$ and $t>0$ and there are constants $\sigma,\mu>0$ and $\lambda<1$ such that 
\begin{equation}\label{hip}
k\bra{t,x,x}< \sigma\frac{\sx^{-\mu}}{t^\lambda}
\end{equation}
then we have the following estimates:\\
Whenever  $\sx^\gamma t > 1$ we have
\begin{equation}
\sigma\frac{\sx^{-\mu}}{t^\lambda}\exp\bra{-c\sx^{{}^{\frac{\bra{1-\frac{N}{2m}}\gamma-\mu}{1-\lambda}}}t}<k\bra{t,x,x}
\end{equation}
and when $\sx^\gamma t < 1$ we have
\begin{equation}
\sigma \frac{\sx^{-\mu}}{t^\lambda}\exp\bra{-c\sx^{{}^{\frac{-\mu}{1-\lambda}}}\,t^{{}^{\frac{\frac{N}{2m}-\lambda}{1-\lambda}}}}
<k\bra{t,x,x}
\end{equation}
where $c$ is some positive constant..
\end{theorem}

We define the following functions for convenience.
$$u\bra{t,x}:=\sigma\frac{\sx^{-\mu}}{t^\lambda}$$
$$ \delta \bra{t,x}:= \frac{k\bra{t,x,x}}{u\bra{t,x}}$$  
We summarise some of the key aspects of Davies' theory \cite{Da4}.
\begin{lemma}[Davies \cite{Da4}] \label{dav}
Given  $0<\alpha,s<1$ and $p$ defined by
\begin{equation}
p+\bra{1-p}\alpha s=s
\end{equation}
$$k\bra{ts,x,x}<\bra{u\bra{\alpha ts,x}}^{1-p}
\bra{u\bra{t,x}}^p\bra{\frac{k\bra{t,x,x}}{u\bra{t,x}}}^p$$
\end{lemma}
\begin{proof}
The three line lemma is applied to the analytic function $z\rightarrow k\bra{z,x,y}$
\begin{eqnarray*}
k\bra{ts,x,x}&<&\bra{k\bra{\alpha ts,x,x}}^{1-p}k\bra{t,x,x}^p\\
&<&\bra{u\bra{\alpha ts,x}}^{1-p}k\bra{t,x,x}^p\\
&=&\bra{u\bra{\alpha ts,x}}^{1-p}\bra{u\bra{t,x}}^p\bra{\frac{k\bra{t,x,x}}{u\bra{t,x}}}^p
\end{eqnarray*}and the proof is complete.
\end{proof}
\begin{definition}
If $\omega_x$ is the distribution such that 
$\ip{f}{\omega_x}=f\bra{x}$ then  we define $G_t\bra{x,x}$
$$G_t\bra{x,x}:= \langle \bra{tH+1}^{-1}\omega_x,\omega_x \rangle = \langle \int\limits_0^{\infty}e^{-\bra{tH+1}s} ds 
\ \omega_x,\omega_x\rangle$$
\end{definition}
\begin{lemma}[Davies \cite{Da4}] \label{dav2}
 \begin{equation}
\frac{G_t\bra{x,x}}{u\bra{t,x}}<\int\limits_0^1 \frac{u\bra{\alpha ts,x}}{u\bra{t,x}}\,\delta \bra {t,x}^p\,ds+e^{-1}\delta\bra {t,x}
\end{equation}
\end{lemma}
\begin{proof}
By lemma \ref{dav} and using the fact that $u\bra{t,x}$ is a decreasing function we have
\begin{equation}
k\bra{ts,x,x}< u\bra{\alpha ts,x}\bra{\frac{k\bra{t,x,x}}{u\bra{t,x}}}^p
\end{equation}
Hence 
$$
 \langle \int\limits_0^{1}e^{-\bra{tH+1}s} ds \ \omega_x,\omega_x\rangle \leq \int\limits_0^{1}u\bra{\alpha ts,x}\bra{\frac{k\bra{t,x,x}}{u\bra{t,x}}}^p ds
$$
Then the estimate 
$$
 \langle \int\limits_1^{\infty}e^{-\bra{tH+1}s} ds \ \omega_x,\omega_x\rangle \leq e^{-1}k\bra{t,x,x}
$$
and dividing by $u\bra{t,x}$ completes the proof.
\end{proof}
\begin{lemma}
$$
G_t\bra{x,x}  = \Sup{g \in Dom{Q}}\{\frac{|g\bra{x}|^2}{t\Quad{g}+\norm{g}^2_2} \, :g \neq 0\}\\
$$
\end{lemma}
\begin{proof}
\begin{eqnarray*}
G_t\bra{x,x} & = & \ip{\bra{tH+1}^{-1}\omega_x}{\omega_x}\\
        & = & \ip {\bra{tH+1}^{-\frac{1}{2}}\omega_x}{\bra{tH+1}^{-\frac{1}{2}}\omega_x}\\
        & = & \Sup{f \in L^2\bra{\Ome}}\{|\ip {\bra{tH+1}^{-\frac{1}{2}}\omega_x}{f}|^2\, : \norm{f}_2=1\}\\
        & = & \Sup{f \in L^2\bra{\Ome}}\{|\bra{tH+1}^{-\frac{1}{2}}f\bra{x}|^2 : \norm{f}_2=1\}\\
        & = & \Sup{g \in Dom{Q}}\{\frac{|g\bra{x}|^2}{\norm{\bra{tH+1}^{\frac{1}{2}}g}^2} \, :g \neq 0\}\\
        & = & \Sup{g \in Dom{Q}}\{\frac{|g\bra{x}|^2}{t\Quad{g}+\norm{g}^2_2} \, :g \neq 0\}\\
\end{eqnarray*}
\end{proof}

\section{Some Inequalities}

\begin{lemma}\label{lem:gamma}
If $\delta$ and $\lambda$ are both positive and less than 1 then
$$\bra{1-\lambda}^{1-\lambda}e^{\lambda-1}\bra{\ln\delta^{-1}}^{\lambda-1}>\delta$$
\end{lemma}
\begin{proof}
Follows from maximising the function $\delta\bra{\ln\delta^{-1}}^{1-\lambda}$ over lambda.
\end{proof}

\begin{lemma} \label{about}
For positive $\alpha<\tfrac{1}{2}$ and $p$ defined as in lemma \ref{dav}
$$\int\limits_0^1 \frac{u\bra{\alpha ts,x}}{u\bra{t,x}}
\,{\delta \bra {t,x}}^p ds+e^{-1}\delta < \bra{\frac{\Gamma\bra{1-\lambda}}{\alpha}+
\bra{1-\lambda}^{1-\lambda}e^{\lambda-2}}
\bra{\ln\delta^{-1}}^{\lambda-1}  $$
\end{lemma}
\begin{proof}
\begin{eqnarray*}
\int\limits_0^1 \frac{u\bra{\alpha ts,x}}{u\bra{t,x}}
\bra{\delta \bra {t,x}}^p ds+e^{-1}\delta
&=&\int\limits_0^1 \bra{\alpha s}^{-\lambda}
\bra{\delta \bra {t,x}}^p ds+e^{-1}\delta
\end{eqnarray*}
substitution of $\tau=\alpha s$ and with $p>\alpha$ we have 

$$\int\limits_0^1 \bra{\alpha s}^{-\lambda}
\bra{\delta \bra {t,x}}^p ds<\frac{1}{\alpha}\int\limits_0^\alpha
\tau^{-\lambda}\exp\bra{-|\ln\delta \,| \tau}d\tau $$
a further change of variable $y=\tau\ln\delta^{-1}$ the RHS becomes
$$ \frac{\bra{\ln\delta^{-1}}^{\lambda-1}}{\alpha}
\int\limits_0^\alpha
y^{-\lambda}e^{-y}dy$$
we use 
$$ \Gamma\bra{1-\lambda}>\int\limits_0^\alpha
y^{-\lambda}e^{-y}dy$$  along with lemma \ref{lem:gamma} to complete the proof. 
\end{proof}
\section{Estimation of the Greens function}
If $g_x\bra{y}$ is an appropriate test function then we can estimate the Greens function with the evaluation
$$\frac{|g\bra{x}|^2}{t\Quad{g}+\norm{g}^2_2}$$ 
We are primarily concerned with 
behaviour of the heat kernel away from the origin, and thus will
consider test functions to reflect this. The choice of test functions is the only significant difference between the analysis 
to find lower bounds for elliptic operators in \cite{Da4} and our analysis here for \Schrodinger operators. We require test functions
that are measurably sensitive to the potential as well as the differential operator. 
\newline
\begin{definition}
 We let $\psi \in C^\infty_c\bra{\R}$ such that $\psi\bra{s}=1$
when $|s|\leq 1$ and $\psi\bra{s}=0$ when $|s|\ge 2$
\end{definition} 

\begin{definition} 
We define $g_x\bra{y}\in \kik{\R^N}$ 
\begin{equation}
g_x\bra{y}=\psi\bra{\frac{x-y}{\sx^{\beta}}}
\end{equation}
for some $\beta\bra{x}$ such that $\beta\bra{x}<1$ for all $x \in \R^N$
\end{definition}
It is straight forward to show that 
$$\norm{\bra{\bra{-\lap}^m g_x}}_\infty\leq c \,\sx^{-\beta 2m}$$
for some poistive  $c$.

\begin{lemma}
There is some positive constant $c>0$ such that 
\begin{equation}\label{greenineq}
 t\Quad{g_x}+\norm{g_x}_2^2  \leq c\sx^{\beta N}\bra{t\sx^{-2m\beta}+t\sx^{\gamma}+1}
\end{equation}
 \end{lemma}
\begin{proof}
There is some positive constant $c$ such that 
$$Q_0\bra{g_x}\le c \sx^{\beta\bra{N-2m}}$$
and moreover there is a positive constant $c$ such that 
$$\norm{V^\frac{1}{2}g_x}^2_2
\leq c\sx^{\gamma +\beta N}$$
\end{proof}

\section {Heat kernel for $t\sx^\gamma>1$}
\begin{lemma}
There is a positive constant $c$ such that 
$$t\Quad{g_x}+\norm{g_x}^2_2\leq c\sx^{{}^{\bra{1-\frac{N}{2m}}\gamma t}}$$
\end{lemma}becomes
\begin{proof}
For $t\sx^\gamma>1$ inequality \eqref{greenineq} reduces to  
\begin{equation}\label{eqn:almost}
t\Quad{g_x}+\norm{g_x}_2^2\leq c\sx^{\beta N}t\bra{\sx^{-2m\beta}+\sx^{\gamma}}
\end{equation}
and we optimise the RHS over $\beta$ with 
$$\beta_M\bra{x}=\frac{1}{2m}\bra{\gamma - \frac{\ln\bra{\frac{2m}{N}-1}}{\ln\bra{\sx}}}$$
and since the definition of $\sx$ implies that
$\frac{2m}{N}-1<\sx^{2m+\gamma}$ it follows that $\beta_M\bra{x}<1$ for all $x$. 
The proof is completed by substituting $\beta_M$ into \ref{eqn:almost} 
\end{proof}

\begin{lemma}\label{later}
There is a positive constant $c$
$$\sigma\frac{\sx^{-\mu}}{t^\lambda}\exp\bra{-c\sx^
{{}^{\frac{\bra{1-\frac{N}{2m}}\gamma-\mu}{1-\lambda}}t}}<k\bra{t,x,x}$$
\end{lemma}
\begin{proof}
From the definition of the test functions we recall that $g_x\bra{x}=1$, hence
\begin{equation}\label{eqn:gh1}
\frac{|g_x\bra{x}|}{t\Quad{g_x}+\norm{g_x}^2_2}> \frac{c}{t}\sx^{{}^{\bra{\frac{N}{2m}-1}\gamma}}
\end{equation}
for some constant $c$.
Application of lemma \ref{about} yields the inequality 

$$c\frac{\sx^{{}^{\bra{\frac{N}{2m}-1}\gamma+\mu}}}{t^{1-\lambda}}<
\bra{\frac{\Gamma\bra{1-\lambda}}{\alpha}+\bra{1-\lambda}^{1-\lambda}e^{\lambda-2}}
\bra{\ln\delta^{-1}}^{\lambda-1}$$
and  simplyfying we have 
$$-c\sx^{\tfrac{\bra{1-\frac{N}{2m}}\gamma-\mu}{1-\lambda}}t
<\ln\delta$$
and our inequality is obtained on inverting the logarithm and substitution of $u\bra{t,x}$.
\end{proof}

\section{Heat kernel for $t\sx^\gamma<1$}
\begin{lemma}
There is some constant $c>0$ such that 
\begin{equation*}
\Quad{g_x}+\norm{g_x}_2^2<ct^{\frac{N}{2m}}
\end{equation*}
\end{lemma}
\begin{proof}
When $t\sx^\gamma>1$ inequality \eqref{greenineq} reduces to  
\begin{equation}\label{eqn:almost2}
\Quad{g_x}+\norm{g_x}_2^2\leq c\sx^{\beta N}\bra{t\sx^{-2m\beta}+1}
\end{equation}
for some $c>0$.
We find that the RHS is optimised over $\beta$ for 
$$\beta_M=\frac{1}{2m\,\ln\sx}\bra{\ln\bra{\frac{2m}{N}-1} + \ln\,t}$$
\end{proof}
\begin{lemma} \label{early}
There is a positive constant $c$ such that 
$$\sigma \frac{\sx^{-\mu}}{t^\lambda}\exp\bra{-c\sx^{{}^\frac{-\mu}{1-\lambda}}\,t^{{}^\frac{\frac{N}{2m}-\lambda}{1-\lambda}}}
<k\bra{t,x,x}
$$
\end{lemma}
\begin{proof}
Using $g_x\bra{x}=1$ we have
\begin{equation}\label{eqn:gh1}
\frac{|g_x\bra{x}|}{t\Quad{g_x}+\norm{g_x}^2_2}> c\,t^\frac{-N}{2m}
\end{equation}
for some constant $c>0$.
Substituting for $u\bra{t,x}$ and an application of lemma \ref{about} we have 
$$c\,\sigma^{-1}\,\sx^{\mu}\,t^{{}^{\lambda-\frac{N}{2m}}}<
\bra{\frac{\Gamma\bra{1-\lambda}}{\alpha}+\bra{1-\lambda}^{1-\lambda}e^{\lambda-2}}
\bra{\ln\delta^{-1}}^{\lambda-1}$$
and we solve to complete the proof.
\end{proof}

\begin{proof}[Proof of theorem \ref{mian}]
 Follows from lemmas \ref{early} and \ref{later}
\end{proof}

\section*{Acknowledgements}
This research was funded by an EPSRC Ph.D grant. I would like to thank Brian Davies
for giving me this problem and his encouragement since. I thank Owen Nicholas and Mark Owen
 for their invaluable discussions. I am indebted 
 to Anita for all her support.

%%%%%%%%%%%%%%%%%%%%%%%%%%%%%%%%%%%%%%%%
\vskip 0.3in
%Department of Mathematics \newline
%Strand \newlines \newline
%Strand \newline
%London WC2R 2LS \newline
%King's College \newline
%England \\
%e-mail: nclaire@mth.kcl.ac.uk
%\vfill
\end{document}